\documentclass[12pt]{amsart}
\usepackage{amssymb}
\usepackage{amsmath}
\usepackage{indentfirst,latexsym}
\usepackage{mathrsfs}
\usepackage{graphicx}
\usepackage{fancyhdr}
\usepackage[usenames]{color}
\date{}
\begin{document}

\title[]{The Eulerian numbers on restricted centrosymmetric permutations}
\author{Marilena Barnabei, Flavio Bonetti, and Matteo Silimbani}
\address{
 Department of Mathematics, University of Bologna\\ P.zza di Porta San Donato 5, 40126 Bologna}
 \email{barnabei@dm.unibo.it\\bonetti@dm.unibo.it\\silimban@dm.unibo.it}
 \maketitle
\noindent {\bf Abstract.}  We study the descent distribution over
the set of centrosymmetric permutations that avoid the pattern of
length  $3$. Our main tool in the most puzzling case, namely,
$\tau=123$ and $n$ even, is a bijection that associates a Dyck
prefix of length $2n$ to every centrosymmetric permutation in
$S_{2n}$ that avoids $123$.
\newline

\noindent {\bf Keywords:} restricted permutations, centrosymmetric
permutations, Dyck prefixes, Eulerian numbers.
\newline

\noindent {\bf AMS classification:} 05A05, 05A15, 05A19.


\section{Introduction}

\noindent A permutation $\sigma\in S_n$ is \emph{centrosymmetric}
if $\sigma(i)+\sigma(n+1-i)=n+1$ for every $i=1,\ldots,n$.
Equivalently, $\sigma$ is centrosymmetric whenever
$\sigma^{rev}=\sigma^c$, where $rev$ and $c$ are the usual reverse
and complement operations. The subset $C_n$ of centrosymmetric
permutations is indeed a subgroup of $S_n$ that, in the even case,
is isomorphic to the hyperoctahedral group $B_{\frac{n}{2}}$, the
natural $B$-analogue of the symmetric group.

\noindent Centrosymmetric permutations have been extensively
studied in recent years from different points of view. For
example, the present authors \cite{bsb} studied the descent
distribution (or Eulerian distribution) over the subset of
centrosymmetric involutions, while Guibert and Pergola \cite{gp}
and Egge \cite{eggedue} studied some properties of $C_n$ from the
pattern avoidance perspective.

\noindent In this paper we merge these two points of view, and
analyze the descent distribution over the set $C_n(\tau)$ of
centrosymmetric permutations that avoid a given pattern $\tau\in
S_3$.

\noindent As well known, the six patterns in $S_3$ are related as
follows:
\begin{itemize}
\item $321=123^{rev}$,
\item $231=132^{rev}$,
\item $213=132^c$,
\item $312=(132^{c})^{rev}$.
\end{itemize}
Since a permutation $\sigma$ is centrosymmetric whenever
$\sigma^{rev}$ and $\sigma^c$ are centrosymmetric, in order to
determine the distribution of the descent statistic over
$C_n(\tau)$, for every $\tau\in S_3$, it is sufficient to examine
the distribution of descents over the two sets
$C_n(132)$ and $C_n(123)$.\\

\noindent In both cases, our starting point is the
characterization of the elements in $C_n(\tau)$, already appearing
in \cite{eggedue}. In the case $\tau=132$, this characterization
allows us to easily determine the descent distribution.\\

 \noindent The case $\tau=123$ presents some more challenging aspects.
 First of all, we observe that the sets
$C_{2k}(123)$ and $C_{2k+1}(123)$ have substantially different
features. In fact, the set $C_{2k+1}(123)$ is in bijection with
the set $S_k(123)$ of $123$-avoiding permutations. In this case,
the descent distribution over $C_{2k+1}$ can be trivially deduced
from the descent distribution over $S_k(123)$, appearing in
\cite{sbb}.\newline

\noindent In the even case, we define a bijection $\Phi$ between
the set of centrosymmetric permutations in $C_{2n}(123)$ and the
set of Dyck prefixes of length $2n$. The map $\Phi$ yields a
bijective proof of the result $|C_{2n}(123)|={2n\choose n}$, that
have been proved in \cite{eggedue} with enumerative techniques.

\noindent Moreover, the bijection $\Phi$ reveals to be a powerful
tool in determining the descent distribution over $C_{2n}(123)$.
In fact, the Dyck prefix $\Phi(\sigma)$ can be split according to
its last return decomposition into subpaths that are either Dyck
paths or elevated Dyck prefixes, namely, Dyck prefixes with no
intersections with the $x$-axis, apart from the origin. The study
of the descent distribution over the sets of permutations that
correspond to Dyck prefixes of these two kinds leads to an
explicit expression of the bivariate generating function
$$T(x,y)=\sum_{n\geq 0}\sum_{\sigma\in C_{2n}(123)}
x^ny^{\textrm{des}(\sigma)}.$$

\section{Preliminaries}

\subsection{Permutations}

\noindent Let $\sigma\in S_n$ and $\tau\in S_k$, $k\leq n$, be two
permutations. We say that $\sigma$ \emph{contains} the pattern
$\tau$ if there exists a subsequence
$\sigma(i_1)\,\sigma(i_2),\ldots\sigma(i_k)$, with $1\leq
i_1<i_2<\cdots<i_k\leq n$, that is order-isomorphic to $\tau$. We
say that $\sigma$ \emph{avoids} $\tau$ if $\sigma$ does not
contain $\tau$. Denote by $S_n(\tau)$ (respectively $C_n(\tau)$)
the set of $\tau$-avoiding
permutations in $S_n$ (resp. $C_n$), where $C_n$ denotes the set of centrosymmetric permutations in $S_n$.\\

\noindent We recall that, given a permutation $\sigma\in S_n$, one
can partition the set $\{1,2,\ldots,n\}$ into intervals
$I_1,\ldots,I_t$, with $I_j=\{k_j,k_j+1,\ldots,k_j+h_j\}$,
$h_j\geq 0$, such that $\sigma(I_j)=I_j$ for every $j$. The
restrictions of $\sigma $ to the intervals in the finest of these
decompositions are called the \emph{connected components} of
$\sigma$. A permutation $\sigma$ with a single connected component
is called \emph{connected}. A permutation is called \emph{right
connected} if $\sigma^{rev}$ is connected.
 The notion of
right connected component of a permutation is defined in the
obvious way.

\noindent For example, the permutation $$\rho=2\ 7\ 6\ 1\ 3\ 5\
4$$ is right connected, while $$\sigma=5\ 7\ 6\ 4\ 2\ 1\ 3$$ is
not.\\

\noindent Note that the right connected components of a
centrosymmetric permutation are mirror symmetric.

\noindent In the following example, the centrosymmetric
permutation $\tau$ is split into its connected components:
$$\tau=7\ 8\ |6\ |4\ 5\ |3\ |1\ 2.$$

\noindent We say that a permutation $\sigma$ has a \emph{descent}
at position $i$ if $\sigma(i)>\sigma(i+1)$. The set of descents of
$\sigma$ is denoted by Des$(\sigma)$, while des$(\sigma)$
indicates the cardinality of Des$(\sigma)$.

\noindent Observe that the descent set of a permutation $\sigma\in
C_n$ must be mirror symmetric, namely, $i\in\,$Des$(\sigma)$
whenever $n-i\,\in\,$Des$(\sigma)$.

\subsection{Lattice paths}

\noindent A \emph{Dyck prefix} is a lattice path in the integer
lattice $\mathbb{N}\times\mathbb{N}$ starting from the origin,
consisting of up-steps $U=(1,1)$ and down steps $D=(1,-1)$, and
never passing below the x-axis. It is well known (see e.g.
\cite{gtouf}) that the number of Dyck prefixes of length $n$ is
$\left({n\atop \left\lfloor\frac{n}{2}\right\rfloor}\right)$. A
Dyck prefix ending at ground level is a \emph{Dyck path}. If it is
not the case, it will be called a
\emph{proper} Dyck prefix.\\

\noindent A \emph{return} of a Dyck prefix is a down step ending
on the $x$-axis. Needless to say, a Dyck prefix is a Dyck path
whenever it has a return at the last position. We say that a Dyck
prefix is \emph{elevated} if either it has no return, or it has
only one return at the last position.\\

\noindent We observe that a given a Dyck prefix $\mathscr{D}$ can
be classified according to the position of its last return
(\emph{last return decomposition}). The path $\mathscr{D}$ can
be:\begin{itemize}
\item a Dyck path \item an elevated proper Dyck prefix \item the
juxtaposition of a Dyck path and an elevated proper
prefix.\end{itemize}

\begin{figure}[h]
\begin{center}
\includegraphics[bb=20 555 416 817,width=.8\textwidth]{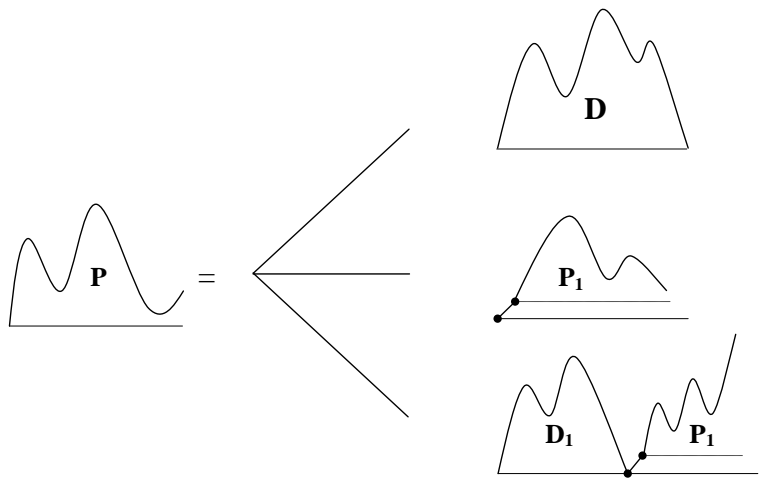} \caption{The last return decomposition of a Dyck prefix.}\label{veniamo}
\end{center}
\end{figure}

\section{The descent distribution over the set $C_n(132)$}

\noindent We begin with two straightforward considerations about
centrosymmetric permutations avoiding $132$.
\begin{itemize}
\item a permutation $\sigma$ belongs to $C_n(132)$ if and only if
the sequence $\sigma(1)\ldots\sigma(n)$ is either $1\,2\ldots n$,
or a sequence of the following kind
$$y\, y+1\,  \ldots\, n\, \beta\, 1\, 2\,  \ldots\, n+1-y$$
where $y>\left\lceil\frac{n}{2}\right\rceil$ and $\beta$, after
renormalization, is either empty or a permutation in
$C_{2y-2-n}(132)$. For example, the eight permutations in
$C_6(132)$ are
$$\begin{array}{lllll}123456\quad&&\qquad& 456123\quad&(\beta=\emptyset)\\
563412\quad&(\beta=12)&\qquad&564312\quad&(\beta=21)
\\ 623451\quad&(\beta=1234)&\qquad& 645231\quad&(\beta=3412)\\ 653421\quad&(\beta=4231)&\qquad&
654321\quad&(\beta=4321)\end{array}$$
\item the set $C_{2n}(132)$ corresponds bijectively to the set
$C_{2n+1}(132)$. In fact, every permutation $\sigma\in
C_{2n}(132)$ corresponds to the permutation $\alpha\in
C_{2n+1}(132)$ defined as follows:
$$\alpha(i)=\left\{\begin{array}{ll}\sigma(i)&\textrm{if }i\leq
n,\\
n+1&\textrm{if }i=n+1\\\sigma(i-1)&\textrm{if
}i>n+1\end{array}\right.$$

\noindent For example, $C_7(132)$ contains the following eight
permutations:
$$\begin{array}{lll}1234567&\qquad\qquad\qquad& 5674123\\
6734512&\qquad\qquad\qquad&6754312
\\ 7234561&\qquad\qquad\qquad& 7564231\\ 7634521&\qquad\qquad\qquad&
7654321\end{array}$$
\end{itemize}
\noindent Denote by $q_{n,k}$ (respectively $r_{n,k}$) the number
of elements in $C_{2n}(132)$ (resp. $C_{2n+1}(132)$) with $k$
descents, and by $$Q(x,y)=\sum_{n\geq 0}\sum_{\sigma\in
C_{2n}(132)} x^ny^{\textrm{des}(\sigma)}=\sum_{n,d\geq 0}
q_{n,d}x^ny^d,$$
$$R(x,y)=\sum_{n\geq 0}\sum_{\sigma\in C_{2n+1}(132)}
x^ny^{\textrm{des}(\sigma)}=\sum_{n,d\geq 0} r_{n,d}x^ny^d,$$ the
generating functions of the two sequences.\\

\noindent Consider the even case. First of all, $q_{0,0}=1$ and
$q_{n,0}=q_{n,1}=1$ for every $n>0$. Moreover, the above
characterization for the elements in $C_{2n}(132)$ yields the
following recurrence for $q_{n,k}$, with $k\geq 2$:
$$q_{n,k}=\sum_{i=1}^{n-1}q_{i,k-2}.$$
These considerations imply immediately the following:
\newtheorem{inizio}{Theorem}
\begin{inizio} We have:
$$Q(x,y)=\frac{x(1+y)}{1-x(1+y^2)}.$$
Hence, for every $n\geq 1$, $$q_{n,k}={n-1\choose
\left\lfloor\frac{k}{2}\right\rfloor}.$$
\end{inizio}
\begin{flushright}
\vspace{-.5cm}$\diamond$
\end{flushright}

\noindent Now we turn to the odd case. A permutation $\alpha\in
C_{2n+1}(132)$ corresponds to a unique permutation $\sigma\in
C_{2n}(132)$. Observe that, if $\sigma$ has an odd number of
descents, then one of these descents is placed at the middle
position, and hence $\alpha$ has an additional descent. In the
other case, $\sigma$ and $\alpha$ have the same number of
descents. These considerations imply that $r_{n,0}=1$ and
$$r_{n,k}=\left\{\begin{array}{ll}q_{n,k}+q_{n,k-1}&\textrm{if }k\textrm{ is even}\\0&\textrm{if }k\textrm{ is odd}\end{array}\right.$$
for every $k\geq 1$. This yields the following:

\newtheorem{eildue}[inizio]{Theorem}
\begin{eildue} We have:
$$R(x,y)=\frac{x}{1-x(1+y^2)}.$$
Hence, for every $n\geq 1$,
$$r_{n,k}=\left\{\begin{array}{ll}\left({n\atop \frac{k}{2}}\right)&\textrm{if
}k\textrm{ is even}\\0&\textrm{if }k\textrm{ is
odd}\end{array}\right.$$
\end{eildue}
\begin{flushright}
\vspace{-.5cm}$\diamond$
\end{flushright}

\section{Characterization of the set $C_n(123)$}

\noindent The characterization of centrosymmetric $123$-avoiding
permutations on an odd number of objects is quite simple. In fact,
we recall that every permutation $\sigma\in C_{2n+1}$ has a fixed
point at $n+1$. Hence, $\sigma\in C_{2n+1}$ avoids $123$ whenever
it has the following structure:

$$\sigma= \alpha'\ \textrm{\emph{n}+1}\ \,\alpha,$$
where $\alpha$ is an arbitrary $123$-avoiding permutation on
$\{1,2,\ldots,n\}$ and $\alpha'$ is the sequence of the
complements to $2n+2$ of the integers $\alpha(n)\cdots\alpha(1)$.

\noindent For instance, if $\alpha=7\ 6\ 4\ 3\ 2\ 1\ 5$, we have
$\sigma=11\ 15\ 14\ 13\ 12\ 10\ 9\ 8\ 7\ 6\ 4\ 3\ 2\ 1\ 5$.\\

\noindent Denote by $v_{n,k}$ the number of permutations in
$C_{2n+1}(123)$ with $k$ descents and by
$$V(x,y)=\sum_{n\geq 0}\sum_{\sigma\in C_{2n+1}(123)}
x^ny^{\textrm{des}(\sigma)}=\sum_{n,d\geq 0} v_{n,d}x^ny^d,$$ the
bivariate generating function of the sequence $v_{n,k}$.

\newtheorem{altro}[inizio]{Proposition}
\begin{altro}
The series $V(x,y)$ has the following explicit expression:
\begin{equation}
V(x,y)=\frac{-1+\sqrt{1-4xy^2-4x^2y^2+4x^2y^4}}{2xy^2(-1-x+xy^2)}\label{uan}
\end{equation}
\end{altro}

\noindent \emph{Proof.} Previous arguments show that the integer
$v_{2n+1,2k+2}$ equals the number of permutations in $S_n(123)$
with exactly $k$ descents. Hence, we have:
$$V(x,y)=1+y^2(E(x,y)-1),$$
where $E(x,y)$ is the generating function of the Eulerian numbers
over $S_n(123)$. It is shown in \cite{sbb} that
$$E(x,y)=\frac{-1+2xy+2x^2y-2xy^2-4x^2y^2+2x^2y^3+\sqrt{1-4xy-4x^2y+4x^2y^2}}{2xy^2(xy-1-x)}.$$
Trivial computations lead to Identity (\ref{uan}).
\begin{flushright}
\vspace{-.5cm}$\diamond$\vspace{.5cm}
\end{flushright}

\noindent We turn now to the even case, and characterize the
elements of $C_{2n}(123)$ by means of the well known decomposition
of a permutation according to its left-to-right minima (recall
that a permutation $\sigma$ has a \emph{left-to-right minimum} at
position $i$ if $\sigma(i)\leq\sigma(j)$ for every $j\leq i$).

\noindent First of all, we observe that a centrosymmetric
permutation $\sigma\in C_{2n}$ is completely determined by its
first $n$ values, namely, by the word
$$w(\sigma)=\sigma(1)\,\sigma(2)\,\ldots\,\sigma(n),$$ and that $w(\sigma)$
can be written as:
$$w(\sigma) = x_1\, w_1\, x_2\, w_2\,\ldots x_k\, w_k,$$ where the integers
$x_i$ are the left-to-right minima of $\sigma$ appearing within
the first $n$ positions and $w_j$ are (possibly empty) words.
Denote by $l_i$ the length of the word
$w_i$.\\


\noindent In order to characterize the elements of $C_{2n}(123)$,
we define a family of alphabets $A_0, A_1,\ldots$ as
follows:\begin{itemize}
\item $A_0=\{1,2,\ldots,2n\}$, \item $A_i$, with $i>0$, is obtained from $
A_{i-1}$ by removing:\begin{itemize}\item the integer $x_i$ and
its complement $2n+1-x_i$, and \item the integers appearing in
$w_i$ together with the corresponding
complements.\end{itemize}\end{itemize}

\noindent For every set $A_i=\{s_1,s_2,\ldots,s_{2h_i}\}$,
$s_1<s_2<\cdots<s_{2h_i}$, we single out its \emph{middle element}
$m(A_i)=s_{h_i}$.

\noindent We have now immediately the following characterization
of the permutations in $C_{2n}(123)$:

\newtheorem{numerino}[inizio]{Proposition}
\begin{numerino}
A centrosymmetric permutation $\sigma$ avoids $123$ if and only if
$$w(\sigma)=x_1\,w_1\,w(\sigma'),$$
where
\begin{itemize}
\item $x_1\geq n$,
\item $w_1=2n\,2n-1\,\ldots\,2n-l_1+1$, with $2n-l_1+1>x_1$,
\item $\sigma'$ is a
centrosymmetric $123$-avoiding permutation over the alphabet
$A_1$,\item the first entry in $\sigma'$ is less than $x_1$.
\end{itemize}\end{numerino}
\begin{flushright}
\vspace{-.5cm}$\diamond$
\end{flushright}

\noindent As a consequence, we have:

\newtheorem{vabene}[inizio]{Corollary}
\begin{vabene} Let $\sigma$ be a permutation in $C_{2n}(123)$,
with \begin{equation}w(\sigma) = x_1\, w_1\, x_2\, w_2\,\ldots
x_s\, w_s.\label{scompo}\end{equation} Then, for every $i\geq 1$,
\begin{equation}\label{berti}x_i\geq m(A_{i-1})\end{equation}
where $m(A_{i-1})$ is the middle element of the alphabet
$A_{i-1}$.\end{vabene}\noindent  If equality holds in
(\ref{berti}), $x_i$ will be called a \emph{tiny} minimum. It is
easily checked that, if $x_i$ is a tiny
 minimum, the left-to-right minimum $x_j$ is also tiny, for every $j>i$.\\

\noindent For example, consider the permutation $\sigma= 11\ 16\
15\ 9\ 7\ 14\ 13\ 12\ 5\ 4\ 3\ 10\ 8\ 2\ 1\ 6$ in $C_{16}(123)$.
Then:
$$w(\sigma)=\underbrace{11}_{x_1}\ \underbrace{16\
15}_{w_1}\ \underbrace{9}_{x_2}\ \underbrace{7}_{x_3}\
\underbrace{14\ 13\ 12}_{w_3}$$ In this case, $w_2$ is empty, and
$\sigma'$ is the $123$-avoiding permutation, order isomorphic to
$6\ 4\ 10\ 9\ 8\ 3\ 2\ 1\ 7\ 5$,  over the alphabet
$A_1=\{3,4,5,7,8,9,10,12,13,14\}$. Note that $7$ is the only tiny
 minimum in $\sigma$.
\newline

\section{A bijection with Dyck prefixes}

\noindent We recursively define a map $\Phi:C(123)\to
\mathscr{P}$, where $C(123)$ is the set of centrosymmetric
$123$-avoiding permutations of any finite even length and
$\mathscr{P}$ is the set of finite Dyck prefixes of even length.
This map associates a permutation $\sigma\in C_{2n}(123)$ with a
Dyck prefix of length $2n$ as follows: decompose $w(\sigma)$ as in
Identity \ref{scompo}. The word $x_2\,w_2\,\ldots\,x_s\,w_s$,
after renormalization, is the word $w(\sigma')$ of some
permutation $\sigma'\in C_{2n-2l_1-2}$, where $l_i$ is the length
of the word $w_i$. Now set $k=2n+1-x_1$. Then:

\begin{itemize}
\item if $k<n+1$, then
$$\Phi(\sigma)=U^kD^{l_1+1}\bar{\Phi}(\sigma'),$$ where
$\bar{\Phi}(\sigma')$ is the Dyck prefix obtained from
$\Phi(\sigma')$ by deleting the leftmost $k-l_1-1$ steps;
\item if $k=n+1$, namely, $x_1$ is tiny, then
$$\Phi(\sigma)=U^{k+1}D^{l_1}\hat{\Phi}(\sigma''),$$ where
$\hat{\Phi}(\sigma')$ is the Dyck prefix is obtained from
$\Phi(\sigma')$ by deleting the leftmost $k-l_1-2$ steps.
\end{itemize}

\noindent It is easy to check that the word $\Phi(\sigma)$ is a
Dyck prefix.\\

\begin{figure}[ht]
\begin{center}
\includegraphics[bb=115 637 412 761,width=.6\textwidth]{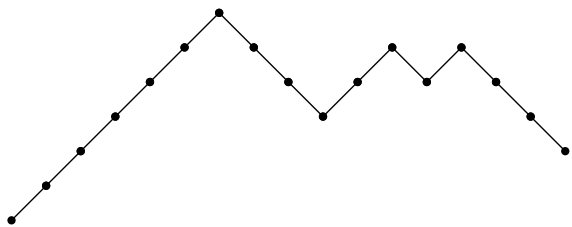} \caption{The
Dyck prefix $\Phi(11\ 16\ 15\ 9\ 7\ 14\ 13\ 12\ 5\ 4\ 3\ 10\ 8\ 2\
1\ 6)$.}\label{info}
\end{center}
\end{figure}

\noindent For example, consider the permutation $\sigma=11\ 16\
15\ 9\ 7\ 14\ 13\ 12\ 5\ 4\ 3\ 10\ 8\ 2\ 1\ 6$. Then,
$\Phi(\sigma)=U^6D^3U^2DUD^3$ (see Figure \ref{info}).\newline

\begin{figure}[h]
\begin{center}
\includegraphics[bb=55 648 357 819,width=.7\textwidth]{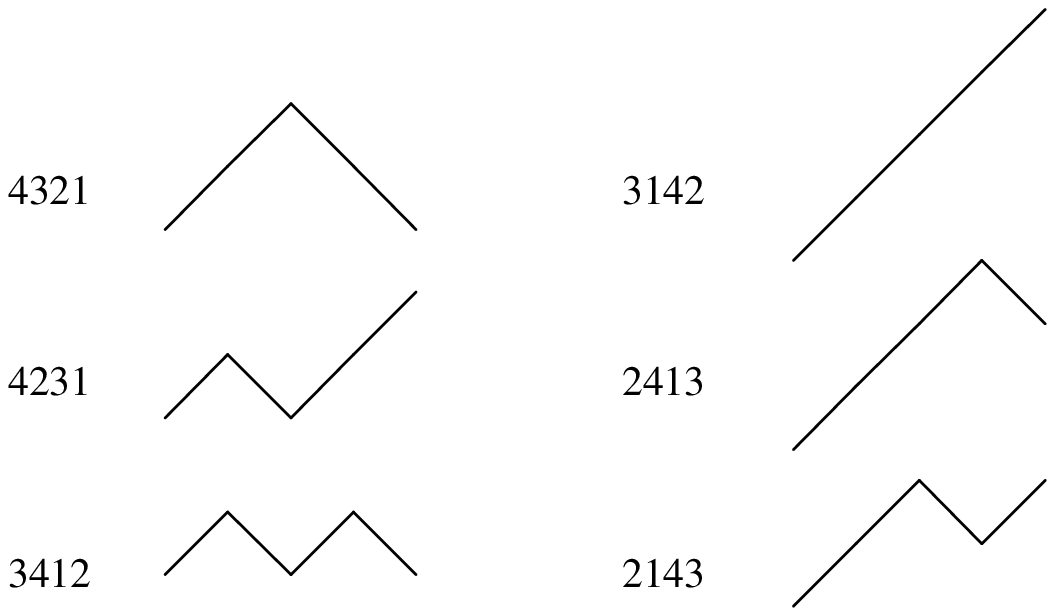} \caption{The Dyck prefixes $\Phi(\sigma)$, with $\sigma\in C_{4}(123)$.}\label{tizcai}
\end{center}
\end{figure}

\noindent In Figure \ref{tizcai}, the prefixes associated with the
six permutations in $C_{4}(123)$ are shown.

\noindent  The map $\Phi$ is a bijection for every positive
integer $n$. In fact, the inverse map $\Phi^{-1}:\mathscr{P}\to
C(123)$ can be recursively defined. Consider a Dyck prefix
$\pi=U^jD^k\pi'$ of length $2n$, where $\pi'$ is a (possibly
empty) lattice path. The permutation $\sigma=\Phi^{-1}(\pi)$ is
defined as follows:

\begin{itemize}
\item if $j\leq n$, set $$\sigma(1)=2n+1-j,\quad \sigma(2)=2n,\quad \sigma(3)=2n-1,\quad\cdots,\quad\sigma(k)=2n-k+2$$
$$\sigma(2n)=j,\quad\sigma(2n-1)=1,\quad\sigma(2n-2)=2,\quad\cdots,\quad\sigma(2n+1-k)=k-1,$$
and let the word $\sigma(k+1)\ldots\sigma(2n-k)$ be the
permutation of the set
$[2n]\setminus\{1,2,\ldots,k-1,j,2n+1-j,2n-k+2,\ldots,2n\}$ that
is order isomorphic to $\Phi^{-1}(U^{j-k}\pi')$;
\item if $j=n+1$, set $$\sigma(1)=n,\quad \sigma(2)=2n,\quad \sigma(3)=2n-1,\quad\cdots,\quad\sigma(k+1)=2n-k+1$$
$$\sigma(2n)=n+1,\quad\sigma(2n-1)=1,\quad\sigma(2n-2)=2,\quad\cdots,\quad\sigma(2n-k)=k,$$
and let the word $\sigma(k+1)\ldots\sigma(2n-k)$ be the
permutation of the set
$[2n]\setminus\{1,2,\ldots,k,n,n+1,2n-k+1,\ldots,2n\}$ that is
order isomorphic to $\Phi^{-1}(U^{j-k-2}\pi')$;
\item if $j>n+1$, set $\sigma(1)=n$, $\sigma(2n)=n+1$, and let $\sigma(2)\ldots\sigma(2n-1)$ be the permutation of
the set $[2n]\setminus\{n,n+1\}$ that is order isomorphic to
$\Phi^{-1}(U^{j-2}D^k\pi')$.
\end{itemize}

\begin{figure}[h]
\begin{center}
\includegraphics[bb=78 658 380 798,width=.5\textwidth]{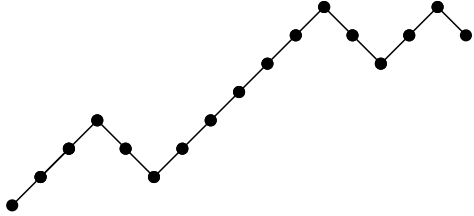} \caption{The
Dyck prefix $U^3D^2U^6D^2U^2D$.}\label{steitare}
\end{center}
\end{figure}

\noindent For example, the permutation associated with the Dyck
prefix $U^3D^2U^6D^2U^2D$ in Figure \ref{steitare} is $\sigma=14\
16\ 8\ 15\ 13\ 7\ 6\ 12\ 5\ 11\ 10\ 4\ 2\ 9\ 1\ 3$.\\

\noindent As an immediate consequence, we obtain the following
result, previously stated in \cite{eggedue}:

\newtheorem{numero}[inizio]{Proposition}
\begin{numero}
The cardinality of the set $C_{2n}(123)$ is the central binomial
coefficient ${2n\choose n }$.\end{numero}
\begin{flushright}
\vspace{-.5cm}$\diamond$
\end{flushright}

\section{Properties of the bijection $\Phi$}

\noindent Some of the properties of a permutation in $C_{2n}(123)$
are related to suitable properties of the associated Dyck prefix.

\noindent First of all, the number of tiny  minima of a
permutation $\sigma$ determines the height of the last point of
its image under $\Phi$. The following result is an immediate
consequence of the definition of the map $\Phi$:

\newtheorem{peradesso}[inizio]{Theorem}
\begin{peradesso}
Let $\sigma$ be a permutation in $C_{2n}(123)$. The $y$-coordinate
of the last point of the path $\Phi(\sigma)$ is twice the number
of tiny  minima in $\sigma$.
\end{peradesso}
\begin{flushright}
\vspace{-.5cm}$\diamond$
\end{flushright}

\noindent In particular, we can characterize the permutations
corresponding to Dyck paths as follows:

\newtheorem{perora}[inizio]{Corollary}
\begin{perora}
Let $\sigma$ be a permutation in $C_{2n}(123)$. The path
$\Phi(\sigma)$ is a Dyck path if and only if $\sigma$ has no tiny
 minimum.
\end{perora}
\begin{flushright}
\vspace{-.5cm}$\diamond$
\end{flushright}

\noindent Observe that the permutation $\sigma$ has no tiny
 minimum if and only if the word $w(\sigma)$ is a
permutation of the set $\{n+1,\ldots,2n\}$. It is easy to verify
that the restriction of $\Phi$ to this set of permutations is a
slightly modified version of the map described
by Krattenthaler in \cite{kratt}.\\

\noindent Consider a permutation $\sigma\in C_{2n}(123)$ with no
tiny  minima. The bijection $\Phi$ is based on a procedure that
associates a Dyck prefix with $\sigma$, processing the word
$w(\sigma)$ from left to right. This allows us to determine at
each step of the procedure the height of the last point of the
lattice path constructed hitherto. More precisely, we will denote
by $P_i$ the Dyck prefix obtained after processing $x_i$ and by
$Q_i$ the Dyck prefix obtained after processing $w_i$. We are
interested in determining the heights $k(P_i)$ and $k(Q_i)$ of the
last point in $P_i$ and $Q_i$, respectively.

\noindent By the definition of $\Phi$, we have $k(P_1)=2n-x_1$ and
$k(Q_1)=2n-x_1-l_1$. Consider now the prefix $P_2$. The alphabet
$A_1$ consists of $2n-2-2l_1$ symbols and $x_2$ is the
$(x_2-1-l_1)$-th smallest element in $A_1$. Hence,
$k(P_2)=(2n-2-2l_1)-(x_2-1-l_1)=2n-1-x_2-l_1$ and
$k(Q_2)=2n-1-x_2-l_1-l_2$. Note that these values do not depend on
$x_1$. By similar arguments, we get the following:
$$k(P_j)=2n-(j-1)-x_j-\sum_{r=1}^{j-1}l_r,$$
$$k(Q_j)=2n-(j-1)-x_j-\sum_{r=1}^{j}l_r.\vspace{1cm}$$

\noindent The previous considerations allow us to relate the right
connected components of permutation $\sigma$ to the returns of the
prefix $\Phi(\sigma)$. More precisely, we have:

\newtheorem{favorevole}[inizio]{Theorem}
\begin{favorevole}
For every $n>0$, the number of right connected components of
$\sigma\in C_{2n}(123)$ is
$$\begin{array}{lll}2\cdot\emph{ret}(\Phi(\sigma))&&\textrm{if }\Phi(\sigma)\textrm{ is a Dyck path}\\ 2\cdot\emph{ret}(\Phi(\sigma))+1&&\textrm{otherwise}\end{array}$$
where \emph{ret}$(\Phi(\sigma))$ is the number of returns of
$\Phi(\sigma)$.
\end{favorevole}

\noindent \emph{Proof} Let $\hat{D}$ be the first return of
$\Phi(\sigma)$, if it exists. Then, if we remove all the steps in
$\Phi(\sigma)$ placed after $\hat{D}$, we obtain a Dyck path
$\hat{\mathscr{D}}$. Such a Dyck path corresponds to a subword
$w'=x_1\,w_1\,\ldots\,x_t\,w_t$ of $w(\sigma)$. By previous
remarks, the integers $x_i$ are non-tiny  minima. Recall that the
last point of $\hat{\mathscr{D}}$ has height
$$k(Q_t)=2n-(t-1)-x_t-\sum_{r=1}^{t}l_r.$$
The path $\hat{\mathscr{D}}$ is a Dyck path whenever $k(Q_t)=0$,
and this is equivalent to
$$x_t=2n-(t-1)-\sum_{r=1}^{t}l_r,$$
which is also equivalent to the fact that the set of the entries
in $w'$ is the interval $[2n+1-z,2n]$, with
$$z=t+\sum_{r=1}^t l_r.$$
Denote by $w(\sigma)=a_1\ldots a_n$. Then, the subwords
$w'=a_1\ldots a_z$ and  $w''=2n+1-a_z\ldots2n+1-a_1$ are connected
components of the permutation $\sigma$.

\noindent Then, we remove from $\sigma$ the two subwords $w'$ and
$w''$, and we obtain a new permutation $\tilde{\sigma}.$ We repeat
this process ret$(\Phi(\sigma))$ times, ending with a Dyck prefix
that is either empty or with no returns. In the first case, the
number of connected components of $\sigma$ is
$2\cdot$ret$(\Phi(\sigma))$. The above considerations imply that
in the second case we get a further connected component.
\begin{flushright}
\vspace{-.5cm}$\diamond$
\end{flushright}

\section{The Eulerian distribution on $C_{2n}(123)$}

\noindent We now study the distribution of the descent statistic
over the set $C_{2n}(123)$. To this aim, we consider the bivariate
generating function
$$T(x,y)=\sum_{n\geq 0}\sum_{\sigma\in C_{2n}(123)} x^ny^{\textrm{des}(\sigma)}=\sum_{n,d\geq 0}
t_{n,d}x^ny^d,$$ where $t_{n,d}$ is the number of permutations in
$C_{2n}(123)$ with $d$ descents.

\noindent Recall that the descent set of a permutation $\sigma\in
C_{2n}(123)$ must be mirror symmetric. This implies that:

$$\textrm{des}(\sigma)=\left\{\begin{array}{lll}2\cdot\textrm{des}(w(\sigma))&&\textrm{if }\sigma(n)\leq n\\
2\cdot\textrm{des}(w(\sigma))+1&&\textrm{otherwise.}\end{array}\right.$$

\noindent The bijection $\Phi$ described and studied in the
previous sections reveals to be an effective tool in the analysis
of the Eulerian distribution on the set  $C_{2n}(123)$. In fact,
it is possible to formulate the
condition that $\sigma$ has a descent at a given position in terms of the associated Dyck path.\\

\noindent We begin with the case of permutations corresponding to
those Dyck prefixes that are the elementary blocks in the last
return decomposition. More precisely:
\begin{itemize}
\item the
set $K_{2n}$ of the permutations in $C_{2n}(123)$ such that
$\Phi(\sigma)$ is a Dyck path. In this case, we denote by
$k_{n,d}$ the corresponding Eulerian number and by $N(x,y)$ the
bivariate generating function
$$N(x,y)=\sum_n\sum_{\sigma\in K_{2n}} x^ny^{\textrm{des}(\sigma)}=\sum_{n,d\geq 0} k_{n,d}x^ny^d,$$
\item the set $CK_{2n}$ of the
permutations in $C_{2n}(123)$ such that $\Phi(\sigma)$ is an
elevated Dyck path. We denote by $ck_{n,d}$ the corresponding
Eulerian number and by $CN(x,y)$ the bivariate generating function
$$ CN(x,y)=\sum_n\sum_{\sigma\in
CK_{2n}} x^ny^{\textrm{des}(\sigma)}=\sum_{n,d\geq 0}
ck_{n,d}x^ny^d,$$
\item the set $G_{2n}$ of the permutations in $C_{2n}(123)$ such that
$\Phi(\sigma)$ is an proper elevated Dyck prefix. We denote by
$g_{n,d}$ the corresponding Eulerian number and by $S(x,y)$ the
bivariate generating function
$$S(x,y)=\sum_n\sum_{\sigma\in
G_{2n}} x^ny^{\textrm{des}(\sigma)}=\sum_{n,d\geq 0}
g_{n,d}x^ny^d.$$
\end{itemize}

\noindent First of all we study the relations between the two
generating functions $N(x,y)$ and $CN(x,y)$. Note that every
permutation $\sigma\in K_{2n}$ has a descent at position $n$.
Moreover:

\newtheorem{fotografia}[inizio]{Proposition}
\begin{fotografia}\label{sfrutt}
Let $\sigma$ be a permutation in $K_{2n}$. The number of descents
of $\sigma$ is
$$\textrm{des}(\sigma)=2(k_1+k_2)+1,$$
where $k_1$ is the number of occurrences of $DDD$ (\emph{triple
falls}) in $\Phi(\sigma)$ and $k_2$ is the number of valleys of
$\Phi(\sigma)$.
\end{fotografia}

\noindent \emph{Proof} Let $w(\sigma)=x_1\,w_1\,\ldots\,x_k\,
w_k$. A descent in $w(\sigma)$ may occur in one of the two
following positions:
\begin{itemize}
\item[1.] between two consecutive symbols $a$ and $b$ in same word
$w_i$. These two symbols correspond to two consecutive down steps
in $\Phi(\sigma)$, that are necessarily preceded by a previous
down step. In fact, if $a$ is not the first symbol in $w_i$, then
$a$ is preceded by a symbol $c$, that also corresponds to a down
step. On the other hand, if $a$ is preceded by $x_i$ in
$w(\sigma)$, then $x_i$ corresponds to the collection of steps
$U^k D$, since $x_i$ can not be tiny, as remarked in the previous
section;
\item[2.] before every left-to-right
minimum $x_i$, except for the first one. These positions
correspond exactly to the valleys of $\Phi(\sigma)$.
\end{itemize}
This implies that des$(w(\sigma))=k_1+k_2$. The assertion now
follows from the previous considerations.
\begin{flushright}
\vspace{-.5cm}$\diamond$
\end{flushright}


\noindent An elevated Dyck path of length $2n$ with $p$ valleys
and $q$ triple falls can be obtained by prepending $U$ and
appending $D$ to a Dyck path of length $2n-2$ of one of the two
following types:
\begin{itemize}
\item[1.] a Dyck path with $p$ valleys and $q$ triple falls, ending
with $UD$,
\item[2.] a Dyck path with $p$ valleys and $q-1$ triple falls, not
ending with $UD$.
\end{itemize}
We note that:
\begin{itemize}
\item[1.] the paths of the first kind are in bijection with Dyck
paths of length $2n-4$ with $p-1$ valleys and $q$ triple falls;
\item[2.] in order to enumerate the paths of the second kind we
have to subtract from the number of Dyck paths of length $2n-2$
with $p$ valleys and $q-1$ triple falls the number of Dyck paths
of semilength $n-1$ with $p$ valleys and $q-1$ triple falls,
ending with $UD$. Dyck paths of this kind are in bijection with
Dyck paths of length $2n-4$ with $p-1$ valleys and $q-1$ triple
falls.\end{itemize} Hence, we have:
\begin{equation}ck_{n,d}=k_{n-1,d-2}-k_{n-2,d-4}+k_{n-2,d-2}\quad (n\geq
2).\label{gorg}\end{equation} In addition, exploiting the last
return decomposition of a Dyck path, we obtain the following
identity, that is a straightforward consequence of Proposition
\ref{sfrutt}:
\begin{equation}\label{addirittura}k_{n,d}=ck_{n,d}+\sum_{i=1}^{n-1}\sum_{j=1}^{d-2}ck_{i,j}k_{n-i,d-1-j}\quad
(n\geq 3),\end{equation}
with the convention $k_{n,d}=0=ck_{n,d}=0$ if $d<0$.\\

\begin{figure}[ht]
\begin{center}
\includegraphics[bb=70 578 459 777,width=.8\textwidth]{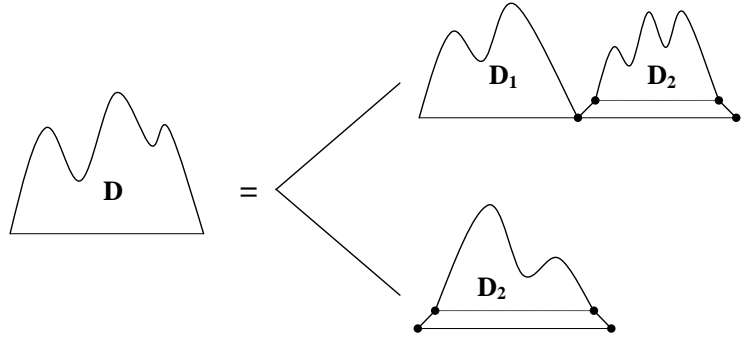} \caption{The last return decomposition of a Dyck path.}\label{andiamo}
\end{center}
\end{figure}

\noindent In fact, as remarked in the previous section, if a Dyck
prefix $\mathscr{D}=\Phi(\sigma)$ is the juxtaposition of a Dyck
path $\mathscr{D}'$ and an elevated proper Dyck prefix
$\mathscr{D}''$, then $w(\sigma)=w'\,w''$, where $w'$ contains the
greatest symbols in $[1,2n]$. Hence,
des$(\sigma)=$des$(\Phi^{-1}(\mathscr{D}'))+$des$(\Phi^{-1}(\mathscr{D}''))+1$.

\noindent Identities (\ref{gorg}) and (\ref{addirittura}) yield:
$$CK(x,y)=xy^2(K(x,y)-1-xy)+x^2y^2(1-y^2)(K(x,y)-1)+1+xy+x^2y,$$
$$K(x,y)=CK(x,y)+y(CK(x,y)-1)(K(x,y)-1).$$
We deduce the following:
$$xy^3(1-xy^2+x)(K(x,y)-1)^2+(2xy^2+2x^2y^2-2x^2y^4-1)(K(x,y)-1)$$$$+xy(1-xy^2+x)=0$$
and hence
\begin{equation}K(x,y)=1+\frac{1-2xy^2-2x^2y^2+2x^2y^4-\sqrt{1-4xy^2-4x^2y^2+4x^2y^4}}{2xy^3(1-xy^2+x)}\label{belli}
\end{equation}
This completes the case of permutations corresponding to Dyck
paths.\\

\noindent Now we turn to the general case. We decompose an
arbitrary Dyck prefix according to its last return, getting

\newtheorem{eunpezzo}[inizio]{Proposition}
\begin{eunpezzo}\label{preocc} For every $n\geq 2$, we have
\begin{equation}t_{n,d}=g_{n,d}+k_{n,d}+\sum_{i=1}^{n-1}\sum_{j\geq
0}g_{i,j}k_{n-i,d-1-j}.\label{certo}\end{equation}
\end{eunpezzo}

\noindent \emph{Proof} If $\sigma$ is neither in $G_{2n}$ nor in
$K_{2n}$, then the Dyck prefix $\Phi(\sigma)$ is the juxtaposition
of a Dyck path $\mathscr{D}'$ and an elevated proper Dyck prefix
$\mathscr{D}''$. In this case, $\sigma$ can be decomposed as:
$$\sigma=\tau_1\,\tau_2\,\tau_3,$$
where the word $\tau_2$, after renormalization, is the permutation
$\alpha=\Phi^{-1}(\mathscr{D}'')$ while $\tau_1\,\tau_3$, after
renormalization, is the permutation
$\beta=\Phi^{-1}(\mathscr{D}')$. This implies that
des$(\sigma)=$des$(\alpha)+$des$(\beta)+2$.
\begin{flushright}
\vspace{-.5cm}$\diamond$
\end{flushright}

\noindent Finally, we express the series $S(x,y)$ in terms of the
functions $T(x,y)$ and $K(x,y)$. Note that, given a Dyck prefix
$\mathscr{D}$ of length $2n-2$, we can prepend to $\mathscr{D}$ an
up step and append either an up or a down step, hence obtaining
two elevated Dyck prefixes $\mathscr{D}'$ and $\mathscr{D}''$.

\begin{figure}[h]
\begin{center}
\includegraphics[bb=91 562 453 793,width=.8\textwidth]{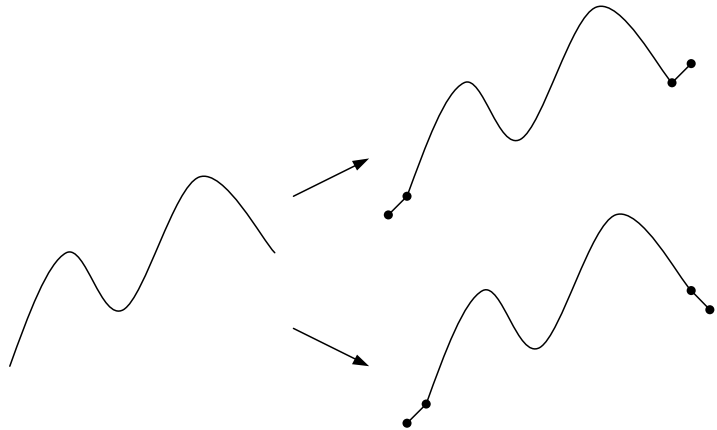} \caption{The generation of two
new Dyck prefixes of length $2n$ from a Dyck prefix of length
$2n-2$.}\label{checapt}
\end{center}
\end{figure}

\noindent The prefix $\mathscr{D}'$ is always proper, while
$\mathscr{D}''$ is proper whenever the prefix $\mathscr{D}$ is not
a Dyck path.

\noindent Denote by $\sigma$ the permutation
$\Phi^{-1}(\mathscr{D})$ and suppose that $\sigma$ has $d$
descents. We want to show that, if we set
$\sigma'=\Phi^{-1}(\mathscr{D}')$ and
$\sigma''\Phi^{-1}(\mathscr{D}'')$, we have:
$$\{\textrm{des}(\sigma'),\textrm{des}(\sigma'')\}=\{d+1,d+2\}.$$

\noindent Note that the number of descents of the permutations
$\sigma'$ and $\sigma''$ depends on the last step in
$\mathscr{D}$:
\begin{itemize} \item if the last step of $\mathscr{D}$ is an up
step, the last entry of the word $w(\sigma)$ is a tiny minimum.
Hence, the word $w(\sigma')$ ends with two consecutive tiny
minima, and des$(\sigma')=$des$(\sigma)+2$. On the other hand,
$w(\sigma'')$ ends with a word $w_k$ of length $1$. Hence,
$\sigma''$ has $d+1$ descents;
\item if the last step of $\mathscr{D}$ is a down step,
in this case, the descent at position $n$ in $\sigma$ splits into
$2$ descents of $\sigma'$. Hence, des$(\sigma')=$des$(\sigma)+1$.
Moreover, neither the last entry of the word $w(\sigma)$ nor the
last entry of the word $w(\sigma'')$ is a left-to-right minimum.
Hence, des$(\sigma'')=$\,des$(\sigma)+2$.
\end{itemize}
Then, we have:
\begin{equation}g_{n,d}=t_{n-1,d-1}+t_{n-1,d-2}-k_{n-1,d-2}\quad (n\geq
2).\label{sino}\end{equation} with the convention $g_{n,d}=0$ and
$t_{n,d}=0$ if $d<0$. Identities (\ref{certo}) and (\ref{sino})
yield the relations:
$$T(x,y)=K(x,y)+S(x,y)-1+y(K(x,y)-1)(S(x,y)-1),$$
$$S(x,y)=1+x+xy(T(x,y)-1)+xy^2T(x,y)-xy^2K(x,y).$$
We deduce the following expression of $T(x,y)$ in terms of
$K(x,y)$:

\newtheorem{etciu}[inizio]{Proposition}
\begin{etciu} We have
\begin{equation}T(x,y)=\frac{-xy^3K^2(x,y)+(1-2xy^2+xy+xy^3)K(x,y)+xy^2-2xy+x}{1-xy+xy^3-xy^2K(x,y)-xy^3K(x,y)}.\label{buoni}
\end{equation}
\end{etciu}
\begin{flushright}
$\diamond$\vspace{.3cm}
\end{flushright}

\noindent An explicit expression for the series $T(x,y)$ can be obtained by combining Identities (\ref{belli}) and (\ref{buoni}).\\

\noindent The first values of the sequence $t_{n,d}$ are shown in
the following table:

$$\begin{array}{l|llllllllll}
n/d&0&1&2&3&4&5&6&7&8&9\\\hline
0&1&&&&&&&&&\\
1&1&1&&&&&&&&\\
2&0&2&3&1&&&&&&\\
3&0&0&3&9&7&1&&&&\\
4&0&0&0&6&20&28&15&1&&\\
5&0&0&0&0&10&50&85&75&31&1\\
\end{array}$$

\end{document}